\documentclass[11pt]{article}
\usepackage{amssymb}
\usepackage{enumerate}
\usepackage{enumitem}
\usepackage{mathrsfs}
\usepackage{latexsym,amsmath,amssymb,amsfonts,epsfig,graphicx,cite,psfrag}
\usepackage{eepic,color,colordvi,amscd}
\usepackage{verbatim}
\usepackage{appendix}

\usepackage{float}
\usepackage{graphicx}
\usepackage{cite}
\usepackage{caption}
\usepackage{subcaption}

\newcommand{\qed}{\hfill $\Box $}
\newcommand{\pf}{\noindent {\bf Proof.} }
\newtheorem{theorem}{Theorem}
\newtheorem{lemma}[theorem]{Lemma}

\newcommand{\ex}{{\rm ex}}
\newcommand{\G}{{\cal G}}
\begin{document}	

	\title{New Bounds on the Anti-Ramsey Number of Independent Triangles\footnote{Research supported by National Key Research and Development Program of China 2023YFA1010200 and  National Natural
Science Foundation of China grant  12271425.}}

\author{
Hongliang Lu, Xinyue Luo and Xinxin Ma \footnote{Corresponding email: maxinxinzzu@163.com}\\School of Mathematics and Statistics\\
Xi'an Jiaotong University\\
Xi'an, Shaanxi 710049, China
}
\date{}


\date{}

\maketitle

		\begin{abstract}
		An edge-colored graph is called \textit{rainbow graph} if all the colors on its edges are distinct. Given a positive integer $n$ and a graph $G$, the \textit{anti-Ramsey number} $ar(n,G)$ is defined to be the minimum number of colors $r$ such that there exists a rainbow copy of $G$ in any exactly $r$-edge-coloring of $K_n$.
		Wu et al. (Anti-Ramsey numbers for vertex-disjoint triangles, \emph{Discrete. Math.}, \textbf{346} (2022), 113123) determined the anti-Ramsey number  $ar(n, kK_3)$ for  $n\geq  2k^2-k+2 $. In this paper, we extend this result by improving the lower bound on $n$ to $n\geq 15k+57$.
		
	\end{abstract}
	\begin{flushleft}
	\end{flushleft}
	\section{Introduction}
	We consider finite graphs without multiple edges or loops. Let  $G$ be a graph with vertex set $V (G)$ and edge set $E(G)$.  We denote the number of edges in $G$ by $e(G)$, i.e., $e(G):=|E(G)|$ and the number of vertices in  $G$ by $|G|$.
	For a subset $X\subseteq V (G)$, we use $G-X$ to denote the graph obtained by deleting all vertices in $X$ and all edges incident to them
	from  $G$. When $X=\{x\}$, we write $G-x$ instead of $G-\{x\}$. Similarly, for $Y\subseteq E(G)$, $G-Y$ denotes the graph obtained by deleting all edges of $Y$ from $G$.
	A subgraph $H$ of $G$ is called  an \emph{induced subgraph} if every pair of vertices in $H$ adjacent in $G$ are also adjacent in $H$. If $H$ and $G$ are graphs such that $G$ does not contain $H$ as a  subgraph, we  say that $G$ is \emph{$H$-free}. Given integers $ n\ge p\ge 1 $, let $ T_p(n) $ denote the \textit{Tur\'an graph}, i.e., the complete $ p $-partite graph on $ n  $ vertices where each part has either $ \lfloor n/p\rfloor $ or $ \lceil n/p \rceil  $ vertices and the edge set consists of all pairs joining distinct parts. Let $t_p(n)$ denote the number of edges in $ T_p(n) $. Let $K_n$ denote the complete graph with $ n $ vertices, and let $[n]$ denote the  set $\{1,2, \cdots ,n \}$. For a given  set $S$ and a positive integer $k$, we use ${S\choose k}$ denote the collection of all possible $k$-element subsets of $S$. Given two positive integers $s$ and $t$, let $K_{s,t}$  denote the complete bipartite graph with partitions of size $s$ and $t$. The notation $G\cong H$ is used to denote that two graphs $G$ and $H$ are isomorphic.
	
	For a vertex $ x\in V(G) $, the neighborhood  of $ x $ in $G$ is denoted by $N_G(x)$. The \textit{degree} of $ x $ in $ G $, denoted by $ d_G(x )$, is the size of $ N_G(x) $. We use $\delta(G)$ and $ \Delta(G) $ to denote the \textit{minimum} and \textit{maximum degrees} of the vertices of $ G $. We use $\bar{d_G}$ to denote the \textit{average degree} of $G$, i.e. $\bar{d_G}=2e(G)/|V(G)|  $.
	
	Graphs in a collection are independent if they share no vertices. Let $X$ and $Y$ be two non-empty subset of $  V(G) $. Let $G[X]$ denote the subgraph \textit{induced} by $X$. Let $E_G(X,Y)$ denote the set of all edges in $G$, with one endpoint in $X$ and one  in $Y$,  and let $e_G(X,Y)=|E_G (X,Y)|$ the number of such edges. When $X=Y$, we simplify the notation to $E_G(X)$ and $e_G(X)$, respectively. For any pair $\{x,y\}\in {V(G)\choose 2}$, define $N_G(x,y)$ as the intersection of the neighborhoods of $x$ and $y$, i.e.,  $N_G(x,y):=N_G(x)\cap N_G(y)$.
	
	For two vertex-disjoint graphs $G,H$, the \textit{join} of $G$ and $H$, denoted by $G \lor H$, is the graph obtained from $G \cup H$ by adding edges joining every vertex of $G$ to every vertex of $H$.
	An \emph{independent vertex set} of $ G $ is a subset of the vertices, no two of which are adjacent.
	
An \emph{$r$-edge-coloring} of a graph assigns $r$ colors to its edges, and an \emph{exactly $r$-edge-coloring} uses all $r$ colors.
	 An edge-colored graph is called \textit{rainbow} if all edges have distinct colors. For a positive integer $ n $ and a graph set $\mathcal{G}$, the \textit{anti-Ramsey number} $ ar(n,\mathcal{G}) $ is the minimum number of colors $r$ such that each edge-coloring of $K_n$ with exactly $r$ colors contains a rainbow copy of $G \in \mathcal{G}$. When $\mathcal{G}=\{G\}$, we denote $ar(n,\mathcal{G})$  by $ ar(n,G) $ for simplicity. The \textit{Tur\'an number} $ \ex(n,G) $ is the maximum number of edges in an $n$-vertex graph with no subgraph isomorphic to $G$.	 A
	graph $H$ on $n$ vertices with $\ex(n, G)$ edges and no copy of $G$ is called an extremal
	graph for $G$. We use $EX(n,G)$ to denote the set of extremal graphs for $G$, i.e.
	\[
	EX(n,G):=\{H\ |\ |V(H)|=n, e(H)= \ex(n,G), \mbox{ and $H$ is $G$-free}\}.
	\]

	The value of $ar(n,G)$ is closely related to  Tur\'an number $\ex(n,G)$ as the following inequality
	\begin{align*}\label{low-upp}
		2+\ex(n,\G) \leq ar(n,G) \leq 1+\ex(n,G),
	\end{align*}
	where $\G=\{G-e\ |\ e\in E(G)\}$.  Erd\H{o}s, Simonovits and S\'os \cite{Erdos1973}  proved there exists a number $n_0(p)$ such that $ar(n,K_p)=t_{p-1}(n)+2$ for $n> n_0(p)$. Montellano-Ballesteros \cite{Mon-Ball2005} and Neumann-Lara \cite{Mon-Ball2002} extended this result to all values of $n$ and $p$ with $n> p\geq 3$.
	
	Let $tK_3$ denote the union of $ t $ independent triangles.
	Let $I_3(G)$ denote the maximum number of independent triangles contained in a given graph $G$. In 1959, Erd\H{o}s and Gallai\cite{erdos1959} determined $ \ex(n,kK_2) $ for any positive integers $n$ and $k$. Later, Erd\H{o}s\cite{erdos1962} proved $ \ex(n,(t+1)K_3)=e(K_t\vee T_2(n-t)) $ for $ n>400t^2 $.
	Moon \cite{Moon1968} proved that $ \ex(n,(t+1)K_3)=e(K_t\vee T_2(n-t))={t\choose 2}+t(n-t)+t_2(n-t)$, and $ K_t\vee T_2(n-t) $ is the unique
	extremal graph that contains no copy of $ (t+1)K_3 $ for $ n>9t/2+4 $. Furthermore,  $$ t_2(n-t)=\lfloor(n-t)^2/4\rfloor=(n-t)^2/4-\epsilon, $$ where
\[\epsilon=\left\{
  \begin{array}{ll}
    0, & \hbox{if $n-t\equiv 0\pmod 2$;} \\
    1/4, & \hbox{otherwise.}
  \end{array}
\right.\]

	\begin{theorem}[Moon, \cite{Moon1968}]\label{extremal}
		If $I_3(G)=t$ and $ n>9/2t+4 $, then
		$$e(G)\leq {t\choose 2}+t(n-t)+\lfloor(n-t)^2/4\rfloor,$$
		with equality holding only if $G\cong K_t\vee K_{\lfloor(n-t)/2\rfloor, \lceil(n-t)/2\rceil}$.
	\end{theorem}

	\begin{lemma}[Moon, \cite{Moon1968}]\label{gamma}
		If $ \nu(G)=h $ and $I_3(G)=0$, then $e(G)\leq h(n-h)$, with equality holding only if $ G\cong K_{h,n-h}$.
	\end{lemma}

	Here we are mainly concerned with the anti-Ramsey problem for
	graphs. The study of anti-Ramsey number for matchings was initiated by Schiermeyer~\cite{Schi2004}, who  proved that \( ar(n, M_s) = ex(n, M_{s-1}) + 2 \) for \( s \geq 2 \) and \( n \geq 3s + 3 \). This result was later improved by Fujita et al. \cite{FKSS}, who extended the range to \( n \geq 2s + 1 \). A complete resolution for all \( s \geq 2 \) and \( n \geq 2s \) was achieved by Chen, Li, and Tu \cite{chen2009}, who determined the exact value of \( ar(n, M_s) \). Interestingly, for the specific case \( n = 2s \), Haas and Young \cite{Haas2012} independently derived the same result using a simpler approach.
Erd\H{o}s, Simonovits and S\'os \cite{Erdos1973} showed that $ ar(n,K_3)=n $. Jahanbekam and West\cite{Jaha2016} determined the anti-Ramsey numbers for $ t $ edge-disjoint perfect matchings in $ K_n $ for $ n \geq 4t+10 $.
	Yuan and Zhang\cite{Yuan} provided the exact results of $ ar(n,tK_3) $ when $ n $ is sufficiently large. Subsequently, Wu and Zhang et al.\cite{Wufangfang} improved this result to $ n\geq 2k^2-k+2 $.
	We aim to establish  the following theorem.	
	\begin{theorem}\label{maintheorem}
		Let $ n,t $ be two integers such that $t\geq 0$ and $ n\geq 15t+57$. The following holds
$$ ar(n,(t+2)K_3)={t\choose 2}+t(n-t)+\lfloor (n-t)^2/4\rfloor+2. $$
	\end{theorem}
	
	\section{Proof of  Theorem \ref{maintheorem}}

	\begin{lemma}\label{pairs}
    Let $ G $ be a $ K_3 $-free graph with $ n $ vertices and  $ e(G)\geq n^2/4-n/2-6t-11 $, where $ n\geq 12t+54 $ and $ t\geq 0 $. Then $G$ contains an independent set $X$  with $|X|\geq n/2-2$ and $ |S|\leq 7 $, where $ S:=\{v\in X\ |\ d_G(v)< \frac{3n}{8}\} $. Furthermore, there exists an independent  set $ X' =X-S$  with $ |X'|\geq n/2-9$ such that for any two vertices $ x,y $ of $X'$,  $ d_G(x)+d_G(y)\geq 3n/4 $.
    \end{lemma}
	\pf Note $ |V(G)|=n $ and $ e(G)\geq n^2/4-n/2-6t-11 $. Then, $ \Delta(G)\geq 2e(G)/n\geq n/2-1-(12t+22)/n\geq n/2-2 $. We choose $ u\in V(G) $ such that $ d_G(u)=\Delta(G) $. Since $ G $ is $ K_3 $-free,  $N_G(u)$ is an independent set of $G$. Let $X:=N_G(u)$. Then $ |X|\geq n/2-2 $.
	
	Suppose that $X$ contains at least eight vertices with $ d_G(v)<3n/8 $. Deleting eight vertices of these vertices yields an induced graph $G'$
 with $ n-8 $ vertices. Then we have
	\begin{align*}
		e(G')> &e(G)- 3n \\
		\geq & n^2/4-n/2-6t-11-3n\\
		=&	\frac{(n-8)^2}{4}+n/2-27-6t\\
		\geq &\frac{(n-8)^2}{4}. 
	\end{align*}
By Theorem \ref{extremal}, $G'$ contains a $K_3$,  contradicting the assumption that $ G $ is $K_3$-free.
	Thus, there are at most seven vertices $v$ in $X$ with $ d_G(v)<3n/8 $, i.e., $ |S|\leq 7 $, where $ S:=\{v\in X\ |\ d_G(v)< \frac{3n}{8}\} $. This completes the proof of the lemma.
	\qed

	\begin{lemma}
		For $ t\geq 0 $	and $ n> 3t+6$.
		$$ar(n,(t+2)K_3)\geq {t\choose 2}+t(n-t)+\lfloor(n-t)/2\rfloor\lceil(n-t)/2\rceil+2.$$	
	\end{lemma}
	\pf Denote $V(K_n)$ by $\{x_1,...,x_n\}$. Consider $ K_t\vee T_2(n-t) $ as a subgraph of $K_n$, where $T_2(n-t)$ the Tur\'an graph on $n-t$ vertices. Define a bijective coloring $ f:E(K_t\vee T_2(n-t))\rightarrow [{t\choose 2}+t(n-t)+\lfloor(n-t)/2\rfloor\lceil(n-t)/2\rceil] $. Let $ G $ be the $ n $-vertex graph $ K_n $ with edge coloring $ f_G $, where
	$$f_G(e)=\begin{cases}
		f(e), & \mbox{if $e \in  E(K_t\vee T_2(n-t))$}, \\
		0,  & \mbox{otherwise}. \\		
	\end{cases}
	$$
	Since $ T_2(n-t) $ contains no rainbow triangles, there are at most $ t $ independent triangles in $ K_t\vee T_2(n-t) $. Thus $ G $ contains no rainbow independent triangles of size $t+2$. Therefore $ ar(n, (t+2)K_3)\geq {t\choose 2}+t(n-t)+\lfloor(n-t)/2\rfloor\lceil(n-t)/2\rceil+2. $ \qed
	\medskip

	\medskip
	\noindent \textbf{Proof of Theorem \ref{maintheorem}.}
	Let 
\[c(n,t)={t\choose 2}+t(n-t)+\lfloor(n-t)/2\rfloor\lceil(n-t)/2\rceil+2.\]
%
Consider a surjective edge-coloring $ c:E(K_n)\rightarrow [c(n,t)] $ of $K_n$ and denote the resulting  edge-colored  $ K_n $ by $H$. Suppose, for contradiction,  that $H$ contains no rainbow $ (t+2)K_3 $. Let $G$ be a rainbow  subgraph of $H$ with exactly $c(n,t)$ edges. For an edge subset $Q\subseteq E(H)$, define $c(Q):=\{c(e)\ |\ e\in Q\}$.
	
	Since $ c(n,t)>ex(n,(t+1)K_3) $,  $G$ contains a copy of
	$ (t+1)K_3 $. Let $M$ denote a set of $ t+1 $ independent triangles of $G$ and let $D$ be the subgraph induced by $ n-3(t+1) $ vertices
that are not part of any triangle in $M$, i.e. $ D=G-V(M)$. The graph  $D$ is $ K_3 $-free,
  otherwise, $G$ would contain  a copy of $ (t+2)K_3 $,  which implies  that $H$ contains $t+2$ independent rainbow triangles (since $G$ is a rainbow graph), a contradiction.  Since $ t\geq 0  $, we may assume that both $V(M)$ and $V(D)$ are non-empty. We say an edge $uv$ is friendly with a vertex $w$ (and vice versa) if $w$ is adjacent to both $u$ and $v$. There cannot be two independent edges in $D$ friendly with two different vertices of a triangle $f$ in $M$, otherwise, $G-V(M\setminus \{f\})$ contains two vertex-disjoint triangles, say $f',f''$, which means that $(M\setminus \{f\}) \cup \{f',f''\}$ is a rainbow $ (t+2)K_3 $ in $ H $, contradicting to the hypothesis. Let $\gamma:=\nu(D)$ and let $\mathcal{M}$ be a maximum matching of $D$.
	
	The triangles of $M$ may be partitioned into two subsets as follows. Define $A$ as the set of triangles $(x,y,z)$ in $M$ where at least one vertex, say $z$, is friendly with at least two independent edges in $\mathcal{M}$. Define $B$ as the set of remaining triangles in $M$.
We denote the number of triangles in $A$ and $B$ by $a$ and $b$, where $ a+b=t+1 $. Then, $M$ can be expressed as $M=\{(x_1,y_1,z_1),...,(x_a,y_a,z_a),(x_{a+1},y_{a+1},z_{a+1}),...,(x_{a+b},y_{a+b},z_{a+b})\}$, where   the first $ a $ triangles belong to $ A $, and the remaining $b$ triangles belong to $ B $.

	If  the triangle $ (x,y,z) $ belongs to $A$ and $z$ is friendly with at least two edges in $\mathcal{M}$, then no vertex in $D$ is friendly with $xy$. Otherwise, suppose that $w\in V(D)$ be friendly with $xy$. Let $x'y'\in \mathcal{M}$ is friendly with $z$ such that $w\notin \{x',y'\}$.
Then  $(M\backslash \{(x,y,z)\})\cup \{(x',y',z),(x,y,w)\}$  is a rainbow $(t+2)K_3$ in $ H $, a contradiction. Therefore, it follows that
\begin{equation}\label{Stru1}
e_G(V(A),V(D))\leq 2a(n-3(t+1)).
\end{equation}
	
	If both  the triangles $ (x_1,y_1,z_1) $ and $ (x_2,y_2,z_2) $ belong to $A$, then neither $x_1$ nor $y_1$ is friendly with $ x_2y_2 $. Otherwise, without loss of generality, assume $x_1$ is friendly with $x_2y_2$. Since there exist two independent triangles of type $ (z_1,p,q) $ and $(z_2,p',q')$, where $pq,p'q'\in \mathcal{M}$,
 then the triangles $ (x_1,y_1,z_1) $ and $ (x_2,y_2,z_2) $ of $M$ could be replaced by the triangles $ (x_1,x_2,y_2) $, $ (z_1,p,q) $ and $(z_2,p',q')$. This replacement forms a set of $t+2$ independent triangles. Consequently, there would be a rainbow  $ (t+2)K_3 $ in $ H $, a contradiction. Therefore, we have
  \begin{align}\label{Stru21}
  e(G[V(A)])\leq {3a \choose 2}-2{a\choose 2},
  \end{align}
and
 \begin{align}\label{Stru2}
  e(G[V(M)])\leq {3(t+1) \choose 2}-2{a\choose 2}.
  \end{align}
	
	Recall that for any $v\in B$, there is at most one  edge in $\mathcal{M}$ friendly with $v$. Therefore, it follows that
\begin{align*}
 e_G(V(B),V(D))&=\sum_{v\in V(B)}e_G(\{v\},V(D))\\
 &\leq3b((n-3(t+1)-2\gamma)+(\gamma+1))\\
 &=3(t+1-a)(n-3(t+1)-\gamma+1),
 \end{align*}
i.e.,
\begin{align}\label{Stru3}
 e_G(V(B),V(D))\leq 3(t+1-a)(n-3(t+1)-\gamma+1).
 \end{align}
Since $I_3(D)=0$, it follows from Lemma \ref{gamma} that
\begin{align}\label{Stru4}
e(D)\leq \gamma(n-3(t+1)-\gamma).
\end{align}
	
	%
	
	\medskip
	\textbf{Claim 1.}~ $a\geq t$.
	
	For the sake of contradiction, suppose that $0\leq a\leq t-1$. Recall that $a+b=t+1$.  Thus we have $b\geq 2$. By  (\ref{Stru1}),  (\ref{Stru2}), (\ref{Stru3}) and  (\ref{Stru4}), we can determine the number of edges in $G$ as follows:
\begin{align*}
		e(G)&=e(G[V(M)])+e_G(V(A),V(D))+e_G(V(B),V(D))+e(D)\\
       &\leq {3(t+1)\choose 2}-2{a\choose 2}+ 2a(n-3(t+1))+3(t+1-a)(n-3(t+1)+1-\gamma)+\gamma(n-3(t+1)-\gamma)\\
		&= 3(t+1)(n-3t/2-1)-a(n+a-3(t+1)+2)+\gamma(n-6(t+1)+3a-\gamma)\\
        &\leq 3(t+1)(n-3t/2-1)-a(n+a-3(t+1)+2)+(n-6(t+1)+3a)^2/4\\
        &=3(t+1)(n-3t/2-1)+5a^2/4+a(n/2-6(t+1)-2)+(n-6(t+1))^2/4.
	\end{align*}
The last expression can be viewed as a function of $a$ that attains its  maximum on the interval  $0\leq a\leq t-1$ when $a=t-1$, provided that  $n\geq 19t/2+35/2$. Therefore,
\begin{align*}
		e(G)\leq & \frac{n^2}{4}+(\frac{t}{2}-\frac{1}{2})n-\frac{t^2}{4}+6t+\frac{61}{4}\\
	<& 	\frac{n^2}{4}+\frac{t}{2}n-\frac{t^2}{4}-\frac{t}{2}+2-\epsilon\quad \mbox{(since $n>13t+57/2+2\epsilon$ and	 $\epsilon\in  \{0,1/4\}) $ }\\
        =&\frac{1}{4}(n-t)^2+t(n-t)+\frac{1}{2}(t^2-t)+2-\epsilon\\
        =&\lfloor\frac{(n-t)^2}{4}\rfloor+t(n-t)+\frac{1}{2}(t^2-t)+2\\
		=& ex(n,(t+1)K_3)+2, 
	\end{align*}
 a contradiction. This completes the proof of Claim 1.
	
	\medskip
	For each triangle $(x_i,y_i,z_i) \in A$, the vertex $z_i$ is friendly with at least two edges of $\mathcal{M}$ denoted by $p_i^1q_i^1$ and $p_i^2q_i^2$. Let $ R:=\cup_{i=1}^{a}\{ p_i^1, q_i^1, p_i^2, q_i^2\} $.
	
	By (\ref{Stru1}),  (\ref{Stru2}) and (\ref{Stru3}), we have
	\begin{align*}
		&e_G(V(M))+e_G(V(M),V(D))
\\= & e_G(V(A),V(D))+e_G(V(B),V(D))+e_G(V(M))\\
		\leq & 2a(n-3(t+1))+{3(t+1)\choose 2}-2{a\choose 2}+3(t+1-a)(n-3(t+1)-\gamma+1),
	\end{align*}
i.e.,
\begin{equation}\label{Eq:1}
		e_G(V(M))+e_G(V(M),V(D))
		\leq 2a(n-3(t+1))+{3(t+1)\choose 2}-2{a\choose 2}+3(t+1-a)(n-3(t+1)-\gamma+1).
\end{equation}	
	By Claim 1, it follows that either $ a=t $ or $ a=t+1 $. Let $n_1:=|V(D)|=n-3(t+1)$.  Next we discuss two cases. 
	
	\medskip	
	\textbf{Case 1.}~ $ a=t+1 $
	
	Then $b=0$ and $A=M$. By (\ref{Eq:1}), we have	
	\begin{align*}
		e_G(V(M))+e_G(V(M),V(D))= & e_G(V(A),V(D))+e_G(V(M))\\
		\leq & 2a(n-3(t+1))+{3(t+1)\choose 2}-2{a\choose 2}\\
		= & 2(t+1)(n-3(t+1))+{3(t+1)\choose 2}-2{t+1\choose 2}\\
		=&2(t+1)n-(5t^2/2+11t/2+3),
	\end{align*}
i.e.,
\begin{align}\label{Ceq:1}
		e_G(V(M))+e_G(V(M),V(D))\leq 2(t+1)n-(5t^2/2+11t/2+3).
	\end{align}
One can see that
	\begin{align*}
		e(D)= &e(G)-e_G(V(M),V(D))-e_G(V(M))\\
		=&c(n,t)-e_G(V(A),V(D))-e_G(A)\quad \mbox{(since $b=0$)}\\
		\geq & \frac{n^2}{4}+\frac{t}{2}n-\frac{t^2}{4}-\frac{t}{2}+2-\epsilon-(2(t+1)n-(5t^2/2+11t/2+3))\quad\mbox{ (by (\ref{Ceq:1}))}\\
		=&	\frac{(n-3(t+1))^2}{4}-(n-3(t+1))/2-t+5/4-\epsilon\\
		\geq& \frac{(n-3(t+1))^2}{4}-(n-3(t+1))/2-t+1,
	\end{align*}
i.e.,
	\begin{align}\label{Ceq:2}
		e(D)\geq\frac{(n-3(t+1))^2}{4}-(n-3(t+1))/2-t+1.
	\end{align}
	
By (\ref{Ceq:2}), we have $ \Delta(D)\geq n_1/2-1-2t/n_1+2/n_1 $. 
	Since $D$ is $K_3$-free, by Lemma \ref{pairs}, there exists an independent set $ X $ in $D$ with $ |X|\geq n_1/2-9 $ such that for any two vertices $ x,y\in X$,  $ d_D(x)+d_D(y)\geq 3n_1/4 $. 	
	 For $\{x,y\}\in {X\choose 2}$, by definition of $N_D(x,y)$, we have $ N_D(x,y)\subseteq V(D)-X $ and  	
	$$|N_D(x,y)|\geq d_D(x)+d_D(y)-(n_1-|X|)\geq -n_1/4+|X|\geq n_1/4-9. $$

	We next consider the edge coloring of $ E(G[X])$ in $ H $.
	
\medskip	
	\textbf{Claim 2.}~For any  $\{u,v\}\in {X\choose 2}$, $ c(uv)\in \{c(e)\ |\ e\in E(M)\} $.
\medskip

Otherwise, suppose that there exist $\{u,v\}\in {X\choose 2}$ and $u'v'\notin  E(M)$ such that $c(uv)=c(u'v')$.
Since the size of the neighborhood $N_D(u,v)$ is at least $n_1/4 - 9$, which is greater than 3, we can select a vertex $w$ from $N_D(u,v)$ where $c(uw)\neq c(u'v')$, and $c(vw)\neq c(u'v')$.
Remember that $G$ is a rainbow graph. From this setup, it's clear that adding the triangle $(u,v,w)$ to $M$ would create $(t+2)$ independent rainbow triangles in $H$, which {\color{blue}{contradicts to}} our earlier assumptions.
This completes the proof of Claim 2.

	Recalled that $ R:=\cup_{i=1}^{a}\{ p_i^1, q_i^1, p_i^2, q_i^2\} $, where  $ (p_i^1,q_i^1,z_i)$ and  $(p_i^2,q_i^2, z_i)$ are two triangles  in $G$ for $ i=1,...,a $.  Since $ |X|\geq n_1/2-9 $, there exists a pair $\{x,y\}\in {X\setminus R\choose 2}$. By Claim 2, there exists a triangle $ (x_j,y_j,z_j) $ in $ A $ such that $ c(xy)\in\{c(e)\ |\ e\in \{x_jy_j,x_jz_j,y_jz_j\}\} $.
 Given that $ |N_D(x,y)|\geq n_1/4-9>3 $, we can select  a vertex $w\in N_D(x,y)\setminus \{p_j^1,q_j^1\}$.  Consequently, the set $ (M-(x_j,y_j,z_j))\cup\{(p_j^1,q_j^1,z_j), (x,y,w)\} $ forms $ (t+2) $ independent rainbow triangles in $H$, which leads to a contradiction.

	\medskip
	\textbf{Case 2.}~ $ a=t$.
	\medskip
	
	Then $b=1$, $ A=\{(x_1,y_1,z_1),...,(x_t,y_t,z_t)\} $ and $ B=\{(x_{t+1},y_{t+1},z_{t+1})\} $. Recalled that $\nu_G(D)=\gamma$.

From equation (\ref{Eq:1}), we derive the following inequality:
	\begin{align*}
		&e_G(V(M))+e_G(V(M),V(D))\\
		\leq&2a(n-3(t+1))+{3(t+1)\choose 2}-2{a\choose 2}+3(t+1-a)(n-3(t+1)-\gamma+1)\\
		=&(2t+3)n-5t^2/2-13t/2-3\gamma-3\quad \mbox{(since $a=t$)},
	\end{align*}
which can be written as:
\begin{align}\label{C2:eq1}
		e_G(V(M))+e_G(V(M),V(D))\leq(2t+3)n-5t^2/2-13t/2-3\gamma-3.
	\end{align}
Using equation (\ref{C2:eq1}), we  observe that
	\begin{align*}
		e(D)= &e(G)-e_G(V(M),V(D))-e_G(V(M))\\
		\geq & \frac{n^2}{4}+\frac{t}{2}n-\frac{t^2}{4}-\frac{t}{2}+2-\epsilon-((2t+3)n-5t^2/2-13t/2-3\gamma-3)\\
		=&\frac{(n-3(t+1))^2}{4}-3(n-3(t+1))/2-3t-7/4-\epsilon+3\gamma\\
		\geq&\frac{(n-3(t+1))^2}{4}-3(n-3(t+1))/2-3t-2+3\gamma,
	\end{align*}
i.e.,
\begin{align}\label{C2:eq2}
		e(D)\geq \frac{(n-3(t+1))^2}{4}-3(n-3(t+1))/2-3t-2+3\gamma.
	\end{align}

	\medskip	
	\textbf{Claim 3.}~$ \gamma\geq n_1/2-(t+3). $
	\medskip	
	
     By contradiction, suppose that $ \gamma\leq n_1/2-(t+4)$.  Since $D$ is $K_3$-free, we can apply Lemma \ref{gamma} to derive the following inequality:
   \begin{equation}\label{C2:eq3}
   e(D)\leq \gamma(n_1-\gamma).
   \end{equation}	
By	combining this result with equations (\ref{C2:eq2}) and   (\ref{C2:eq3}), we conclude that
\begin{align*}
n_1^2/4-3n_1/2-3t-2&\leq \gamma(n_1-\gamma-3)\\
&\leq (n_1/2-t-4)(n_1/2+t+1)\quad \mbox{(since $\gamma\leq n_1/2-t-3$)}\\
&=n_1^2/4-3n_1/2-(t+4)(t+1).
\end{align*}
So we may infer that
\begin{align*}
t^2+2t+2\leq 0,
\end{align*}
which contradicts to that $t\geq 0$.
		This completes the proof of Claim 3.
	\medskip

By equation (\ref{C2:eq2}) and Claim 3, we have
\begin{align}
e(D)\geq \frac{n_1^2}{4}-6t-11.	
\end{align}

   Recall that $D$ is $ K_3 $-free and $n_1=n-3t-3\geq 12t+54$. By Lemma \ref{pairs}, there exists $ X\subseteq V(D) $ such that $X$ forms an independent set in $D$ with $ |X|\geq n_1/2-9 $ and $ d_{D}(x)+d_{D}(y)\geq 3n_1/4 $ for any pair  $ \{x, y\}\in {X\choose 2} $.
	Since for any pair $\{x,y\}\in { X\choose 2} $, $ N_D(x,y)\subseteq V(D)\backslash X $, we can derive
\begin{align}\label{C2:eq4}
|N_D(x,y)|\geq d_{D}(x)+d_{D}(y)-(n_1-|X|)\geq -n_1/4+|X|\geq n_1/4-9\geq 4.
\end{align}

	
%

To complete our proof, we rely on the following  claim.

	\medskip
	\textbf{Claim 4.}~There exists a subset $T\subseteq X$  with $ |T| \geq n_1/2-9-2t $ such that for any pair $\{x,y\}\in {T\choose 2}$,  $ c(xy)\in  \{c(e)\ |\ e\in E(B)\} $.
	\medskip
	
	Consider a triangle  $ (x_i,y_i,z_i) $ from $ A $, we claim that
\[
|\{e\in {X\choose 2}\  |\  c(e)\in \{c(x_iy_i),c(x_iz_i),c(y_iz_i)\}\}|\leq 1.
\]
 Suppose, for contradiction, that there are two distinct pairs $ \{u_1,v_1\},\{u_2,v_2\} \in {X\choose 2}$  such that both $ c(u_1v_1)$ and $c(u_2v_2)$ belong to $c(x_iy_i),c(x_iz_i),c(y_iz_i)\}$.  Since $p_i^1q_i^1,p_i^2q_i^2\in E(D)$ and $X$ is an independent set of graph $D$, there exists some $j\in [2]$ such that
 \[
|\{u_j,v_j\}\cap \{p_i^1,q_i^1,p_i^2,q_i^2\}|\leq 1.
 \]
 Without loss of generality, assume that $\{u_1,v_1\}\cap \{p_i^1,q_i^1\}=\emptyset$.
Given  $ |N_D(u_1, v_1)|\geq n_1/4-9>3 $,
   we can select a vertex $ w \in N_D(u_1,v_1) $ such that  $w\notin \{p_i^1,q_i^1\}$.
   Since $ G $ is a rainbow graph and $ u_1v_1\notin E(D) $, the set $(M\setminus \{(x_i,y_i,z_i)\})\cup \{(p_i^1,q_i^1,z_i),(u_1,v_1,w)\} $ forms  $t+2$ independent rainbow triangles in $ H $, a contradiction.

Let $P:=\{e\in {X\choose 2}\ |\ c(e)\in c(E(A)) \}$. It's clear that $|P|\leq t$. Let $T:=X-V(P)$. Then $T$ is a desired subset,  completing the proof of Claim 4.

	\medskip
	\textbf{Claim 5.}~The set $\{c(e)\ |\ e\in {T\choose 2}\}$ is a singleton subset of $c(E(B))$.
	\medskip
	
	By Claims 2 and 4, we know $\{c(e)\ |\ e\in {T\choose 2}\}\subseteq c(E(B))$. Next  we show that this set contains exactly one color by contradiction. Suppose that $|\{c(e)\ |\ e\in {T\choose 2}\}|\geq 2$.
 Since $ |T|\geq |X|-2t\geq 4 $,  there exist two disjoint pairs $ \{u_1,v_1\}$ and $\{u_2,v_2\}$ in ${T\choose 2}$ such that   $\{c(u_1v_1),c(u_2v_2)\}\in { c(E(B))\choose 2}$. By equation (\ref{C2:eq4}), we may find two different vertices $w_1,w_2$ in $D$ such that $w_i\in N_D(u_i,v_i)$ for $i\in [2]$. This allows us to construct the set $(M\backslash B)\cup \{(u_i,v_i,w_i)\ |\ i\in [2]\}$, which forms $t+2$ vertex-disjoint rainbow triangles in $H$, a contradiction. This completes the proof of Claim 5.

	\medskip
	\textbf{Claim 6.}~For  $v\in \{x_{t+1},y_{t+1},z_{t+1}\} $,  $N_G(v)\cap V(D)$ is an independent set in $G$. 
	\medskip
	
Assume for contradiction that  there exists $v\in  \{x_{t+1},y_{t+1},z_{t+1}\}$ such that $G[N_G(v)\cap V(D)]$ contains an edge $ xy $. Since $ |T|\geq n_1/2-9-2t\geq 4 $, we can choose two vertices $u_1,v_1\in T\backslash \{x,y\}$ and  a vertex $ w \in  N_D(u_1,v_1)\backslash \{x,y\}$. Then the set $(M\backslash B)\cup \{(u_1,v_1,w), (v,x,y)\}$ forms $t+2$ vertex-disjoint triangles in $H$,  a contradiction. This completes the proof of Claim 6.	
	
	\medskip
	
	By Claim 5, we may assume that  for any  $\{x,y\}\in { T\choose 2} $, $ c(xy)=c(x_{(t+1)}y_{(t+1)}) $. Define $G'=G-\{x_{(t+1)}y_{(t+1)}\}$ and $ D'=V(D)\cup\{x_{t+1},y_{t+1},z_{t+1}\}$.
	
\medskip
	\textbf{Claim 7.}~Every  vertex in $V(D)$  is unfriendly with   $ x_{t+1}z_{t+1} $ and $ y_{t+1}z_{t+1} $. 
	\medskip	

Suppose otherwise that there exists a vertex $v'$ which is friendly with $x_{t+1}z_{t+1}$ in $ G $.
Since $ T\geq n_1/4-9-2t\geq 4 $, we can select $\{u,v\}\in {T-v'\choose 2}$ and choose a vertex $ w $ of $ N_D(u,v)\backslash \{v'\}$ (as $ |N_D(u,v)|\geq n_1/4-9\geq 4 $). Since $G$ is a rainbow graph and $ c(uv)=c(x_{(t+1)}y_{(t+1)}) $, the set $ (M\backslash B)\cup \{(u,v,w), (v',x_{(t+1)},z_{(t+1)})\} $ forms  $t+2$ independent rainbow triangles in $ H $, a contradiction. This completes the proof Claim 7.	
	
	By Claims 6 and $7$, $G'[D']$ is $K_3$-free. Applying Lemma \ref{gamma}, we have
\begin{align}\label{bound-e(D')}
 e(G'[D'])\leq \gamma'(n-3t-\gamma'),
 \end{align}
  where $ \gamma'=\nu(G'[D'])$.
		Recall that for $ 1\leq j\leq t $, there is no vertex in $D$, which is friendly with  $x_jy_j$.

 \medskip
	\textbf{Claim 8.}~For $1\leq j\leq t$, every  vertex in $V(B)$  is unfriendly with   $ x_jy_j $.
	\medskip	

 Suppose otherwise that  $ z_{t+1} $ is friendly $x_1y_1$. Since $ |T|\geq n_1/4-9-2t\geq 4 $, we may select  $\{u,v\}\in {T\backslash \{p_1^1,q_1^1\}\choose 2}$  and choose a vertex $ w \in  N_D(u,v)\backslash \{p_1^1,q_1^1\}$  (as $ |N_D(u,v)|\geq n_1/4-9\geq 4 $). Since $ G' $ is a rainbow graph and $ c(uv)=c(x_{(t+1)}y_{(t+1)}) $,  the set  $ (M\backslash (B\cup \{(x_1,y_1,z_1)\})\cup \{(u,v,w), (z_1,p_1^1,q_1^1)\cup (x_1,y_1,z_{t+1}) \}$ forms $ (t+2) $ independent rainbow triangles in $ H $, a contradiction. This completes the proof of Claim 8.
	
	Applying (\ref{Stru1}) and Claim 8,  we may derive
\begin{align*}
 e_{G'}(V(A),D')&=e_{G'}(V(A),V(D))+e_{G'}(V(A),V(B))\\
 &\leq 2t(n-3(t+1))+ 6t=2t(n-3t),
\end{align*}
which can be written as:
\begin{align}\label{C3:eq1}
 e_{G'}(V(A),D')\leq 2t(n-3t).
\end{align}
Thus we have
	\begin{align*}
	  e(G)=&e(G')+1\\
		=&e_{G'}(V(A))+e_{G'}(V(A),D')+e(G'[D'])+1\\
		\leq& {3t\choose 2}-2{t\choose 2}+2t(n-3t)+\gamma'(n-3t-\gamma')+1\quad\mbox{(by (\ref{Stru21}), (\ref{bound-e(D')}) and (\ref{C3:eq1}))}\\
 \leq& {3t\choose 2}-2{t\choose 2}+2t(n-3t)+\lfloor(n-3t)^2/4\rfloor+1\\
		\leq& ex(n,(t+1)K_3)+1<c(n,t),
	\end{align*}
	which contradicts with the fact that $ e(G)= ex(n,(t+1)K_3)+2 $. 	This completes the proof.
	\qed

\end{document}